\documentclass{amsart}

\usepackage{amsmath,amsfonts,amsthm,amssymb,latexsym,amscd,amsopn,eucal,mathrsfs,graphics,color,enumerate,lscape}

\usepackage[colorlinks]{hyperref}

\DeclareMathOperator{\Pic}{Pic}

\DeclareMathOperator{\Br}{Br}

\newcommand{\OO}{\mathcal{O}}
\newcommand{\PP}{\mathbb{P}}
\newcommand{\tns}{\otimes}

\title{Genus One Curves and Brauer-Severi Varieties}
\author{Aise Johan de Jong}
\address{Department of Mathematics, Columbia University, New York, NY 10027}
\email{dejong@math.columbia.edu}
\author{Wei Ho}
\address{Department of Mathematics, Columbia University, New York, NY 10027}
\email{who@math.columbia.edu}

\begin{document}

\maketitle

\section{Introduction}

\noindent
Let $K$ be a field. Let $A$ be a central simple algebra over $K$ and
let $X$ be the associated Brauer-Severi variety over $K$.
It has recently been asked \cite{openproblems} if there exists
a genus $1$ curve $C$ over $K$ such that $K(C)$ splits $A$.
In other words, is there a genus one curve $C$
over $K$ with a morphism $C \to X$? In this short note,
we explicitly construct such a curve in the case where 
$X$ has dimension $\leq 4$ (equivalently, when $A$ has degree $\leq 5$).

\medskip\noindent
There is some related work \cite{KL, ciperiani-krashen}
regarding which Brauer classes are split by a given
(genus $1$) curve over $K$.

\section{Index 2 and Index 3}

\noindent
These cases are covered by previous work \cite{Swets}. We briefly describe 
constructions for these two cases, since the higher cases below are similar in 
spirit.

\medskip\noindent
Let $A$ be a quaternion algebra over $K$, and let $X$ be a genus zero curve 
representing the same Brauer class. Let $L$ be a degree $2$ line bundle on 
$X$, so a section of $L^{\tns 2}$ cuts out a degree $4$ subscheme $D$ of $X$.  
Then there exists a double cover $C$ of $X$ ramified exactly at $D$ by the
branched covering trick. The genus of $C$ is $1$. For a general section,
when the characteristic of $K$ is different from $2$, 
the curve $C$ will be a smooth irreducible genus one curve.

\medskip\noindent
Now let $A$ be a central simple algebra over $K$ of degree $3$, with $X$ the 
corresponding Brauer-Severi variety.  Then the inverse of the canonical bundle 
of $X$ is a line bundle whose general sections cut out genus one curves.

\section{Index 4}

\noindent
Let $A$ be a central simple algebra over $K$ of degree $4$ and let
$X$ be the corresponding Brauer-Severi variety.
Let $\alpha \in \Br(K)$ be the class of $A$, so $\alpha$
is a nontrivial element of index $4$ in  $\Br(K)$ and has period $2$ or $4$.
By \cite[Corollary 15.2.a]{Pierce}, the class of $2\alpha$ has index $2$
or $1$. Let $Y$ be a Brauer-Severi variety of dimension $1$ whose
Brauer class is $2\alpha$.

\medskip\noindent
It is well known that that the intersection of two general sections of 
$\OO_{\PP^3}(2)$ is a smooth irreducible genus one curve. Another way to 
describe this curve is as the zero locus of a general section of the 
pushforward $\pi_* \OO_{\PP^3 \times \PP^1}(2, 1)$, where
$\pi: \PP^3 \times \PP^1 \to \PP^3$ is the first projection.
To generalize this construction for 
our situation, we descend this vector bundle to $X \times Y$.

\medskip\noindent
We claim that the line bundle $\OO_{X_{\bar{K}}}(2) \boxtimes 
\OO_{Y_{\bar{K}}}(1)$ on $X_{\bar{K}} \times Y_{\bar{K}}$ descends to a line 
bundle $\mathcal{L}$ on $X \times Y$.  In other words, we want to show that 
$\OO_{X_{\bar{K}}}(2) \boxtimes \OO_{Y_{\bar{K}}}(1)$ is in the image of the map
$$\Pic(X \times Y) \to \Pic(X_{\bar{K}} \times Y_{\bar{K}}),$$
and more precisely, in the image of the map
$$\Pic(X \times Y) \to \Pic(X_{\bar{K}} \times 
Y_{\bar{K}})^{\mathrm{Gal}(\bar{K}/K)}.$$
The next term in the low degree exact sequence coming from the Leray spectral
sequence is the Brauer group $\Br(K)$.
Similarly, there is an exact sequence
$$\Pic(X) \to \Pic(X_{\bar{K}})^{\mathrm{Gal}(\bar{K}/K)} \to \Br(K).$$
The obstruction to the line bundle $\OO_{X_{\bar{K}}}(1)$ coming from a line 
bundle on $X$ is exactly the class $\alpha$ in $\Br(K)$, so because the 
differential is a homomorphism, the obstruction for $\OO_{X_{\bar{K}}}(2)$ is 
$2\alpha$.  Similarly, the obstruction for $\OO_{Y_{\bar{K}}}(1)$ is $2\alpha$.
By the K\"unneth formula, the obstruction for $\OO_{X_{\bar{K}}}(2) 
\boxtimes \OO_{Y_{\bar{K}}}(1)$ is $2\alpha + 2\alpha = 0$.

\medskip\noindent
Therefore, there exists a line bundle $\mathcal{L}$ on $X \times Y$ as above, 
and the pushforward $\pi_* \mathcal{L}$ via the projection $\pi: X \times Y \to 
X$ is a rank $2$ vector bundle on $X$. By base change, the bundle
$\pi_*\mathcal{L}$ on $X$ has many sections, and a general section
cuts out a genus one curve on $X$.

\section{Index 5}

\noindent
Let $A$ be central simple algebra $A$ over $K$ of degree $5$.
Let $\alpha \in \Br(K)$ be the class of $A$, so
$\alpha$ is a nontrivial element in $\Br(K)[5]$.  Let 
$X$ and $Y$ be Brauer-Severi varieties representing the classes $\alpha$ and 
$2\alpha$, respectively.  We construct a genus one curve in $X$ by finding a 
general section of a vector bundle over $X \times Y$.

\medskip\noindent
The following observation may be found in \cite{dAlmeida} and was explained
to us by Laurent Gruson (private communication). The sheaf
$\OO_{X_{\bar{K}}}(1) \boxtimes \Omega^1_{Y_{\bar{K}}}(2)$
is a rank $4$ vector bundle on
$X_{\bar{K}} \times Y_{\bar{K}}$.
The zero locus of a general section is a closed smooth 
subvariety of $X_{\bar{K}} \times Y_{\bar{K}}$
whose projection to $X_{\bar{K}}$ is a smooth irreducible 
genus one curve. 

\medskip\noindent
We claim that the vector bundle $\OO_{X_{\bar{K}}}(1) \boxtimes 
\Omega^1_{Y_{\bar{K}}}(2)$ on $X_{\bar{K}} \times Y_{\bar{K}}$ descends to a 
vector bundle $\mathcal{E}$ over $K$. Because
$\OO_{X_{\bar{K}}} \boxtimes \Omega^1_{Y_{\bar{K}}}$
certainly descends, we want to show that the line 
bundle $\OO_{X_{\bar{K}} \times Y_{\bar{K}}}(1,2)$
is in the image of the map
$$
\Pic(X \times Y) \to \Pic(X_{\bar{K}} \times 
Y_{\bar{K}})^{\mathrm{Gal}(\bar{K}/K)}.
$$
As in the index $4$ case, the obstruction lies in $\Br(K)$, and an 
almost identical computation gives
$\alpha + 2(2\alpha) = 5\alpha = 0$ in $\Br(K)$.

\medskip\noindent
By base change, the vector bundle $\mathcal{E}$ on
$X \times Y$ has many sections, so we may take a general section as 
above. The projection to $X$ of the zero locus of this section is a genus one curve, 
as desired.

\bibliography{bibg1brauer}
\bibliographystyle{amsplain}

\end{document}